\documentclass[12pt]{article}
\usepackage{amsmath,amsthm,amsfonts,amssymb,MnSymbol,mathrsfs}
\usepackage[english]{babel}

\usepackage[applemac]{inputenc}
\usepackage[all]{xy}
\numberwithin{equation}{section}

  \newcommand{\const}{\rm const}

  \newcommand{\Tail}{\rm  Tail}

\newtheorem{theorem}{Theorem}[section]
\newtheorem{definition}{Definition}[section]
\newtheorem{lemma}{Lemma}[section]

\newtheorem{remark}{Remark}[section]

\newtheorem{example}{Example}[section]
\begin{document}



\title{\textbf{Criterion for the coincidence of strong and weak Orlicz spaces}}

\footnotesize\date{}

\author{Maria Rosaria Formica ${}^{1}$, Eugeny Ostrovsky ${}^2$}

\maketitle

\begin{center}
${}^{1}$ Parthenope University of Naples, via Generale Parisi 13,\\
Palazzo Pacanowsky, 80132,
Napoli, Italy. \\

e-mail: mara.formica@uniparthenope.it \\

\vspace{4mm}

${}^2$ Department of Mathematics and Statistics, Bar-Ilan University, \\
59200, Ramat Gan, Israel. \\

e-mail: eugostrovsky@list.ru\\

\end{center}

\begin{abstract}

We provide necessary and sufficient conditions for the coincidence,
up to equivalence of the norms, between strong and weak Orlicz
spaces. Roughly speaking, this coincidence  holds true only for the
so-called {\it exponential} spaces.\par We find also the exact value
of the embedding constant which appears in the corresponding norm
inequality.
 \end{abstract}

\vspace{2mm}


\vspace{2mm}

 \noindent {\footnotesize {\it Key words and phrases}:

Measure, Orlicz space, Young-Orlicz function, norm equivalence, tail
function and tail norm, expectation, Lorentz spaces,
Orlicz-Luxemburg strong and weak norms, embedding constant, Markov-Tchebychev's inequality.\\

\vspace{2mm}


\vspace{4mm}

\section{Notations. Definitions. Statement of the problem.}

\vspace{4mm}

 Let $ (X = \{x\}, \cal{F}, \mu)$ be a measurable space with atomless sigma-finite non-zero measure $ \mu.$ Let $N = N(u), \ u
 \in \mathbb{R}
 $,
 be a non-negative numerical-valued Young-Orlicz function. This means that $N(u)$ is even, continuous, convex, strictly increasing to infinity as $u\geq
 0$, $u\to \infty$
 and such that
\begin{eqnarray*}
\lim_{u \to 0} \frac {N(u)}{u} = 0 \ , \ \ \ \lim_{u \to \infty}
\frac {N(u)}{u} = +\infty.
\end{eqnarray*}
In particular,
$$
N(u) = 0 \ \Leftrightarrow  \ u = 0.
$$
Denote by $M_0 = M_0(X,\mu) $ the set of all numerical-valued
measurable functions  $f: X \to \mathbb{R}$, finite almost
everywhere.
 The Orlicz space $ L(N) = L(N; X, \mu)$ consists of all functions  $ f: X \to \mathbb{R} $ from the set $ M_0(X,\mu)$
 for which the classical Luxemburg norm $ \ ||f||_{L(N)}$  (equivalent to the Orlicz norm) or, in more
detail, the {\it strong } Luxemburg norm $ ||f||_{sL(N)}$ defined by
\begin{equation}\label{Luxembur norm}
 ||f||_{L(N)}  = ||f||_{sL(N)} := \inf \left\{ k > 0 \, : \, \int_X N(|f(x)|/k) \ d\mu(x) \le 1 \ \right\}
 \end{equation}
is finite.\par

Furthermore, if $0<||f||_{L(N)}<\infty$, then
\begin{equation}\label{integral less-equal 1}
\int_X N \left( \ \frac{|f(x)|}{||f||_{L(N)}} \ \right) \
d\mu(x)\leq 1.
\end{equation}

Note that the equality sign occurs in \eqref{integral less-equal 1}
if in addition the Young - Orlicz function $N(\cdot)$ satisfies the
well known $\Delta_2$-condition. Moreover, if there exists $k_0 >0$
such that $\displaystyle \int_X N \left( \ \frac{|f(x)|}{k_0} \
\right) \ d\mu(x) = 1$, then $f \in L(N)$ and $k_0=||f||_{L(N)}$
 (see \cite{Kras Rut1961}, Chapter 2, Section 9).

\vspace{4mm}

The Orlicz spaces have been extensively investigated by M. M. Rao
and Z. D. Ren in \cite{RaoRenTheory 1991,RaoRenAppl1997}; see also
\cite{Bennett Sharpley1988,Kras
Rut1961,Masta2016,OsMono1999,OsSirHIAT2007}, etc. Recently in
\cite{Fiorenza2018a} (see also \cite{Fiorenza2018b}) the authors
studied the Gagliardo-Nirenberg inequality in rearrangement
invariant Banach function spaces, in particular in Orlicz spaces.
\par

Note that  the so-called {\it exponential} Orlicz spaces are
isomorphic to suitable Grand Lebesgue Spaces, see \cite{Buld Mush
OsPuch1992,KozOs1985,KozOsSir2017,OsMono1999}. For some properties,
variants and applications of the classical Grand Lebesgue Spaces see
for example
\cite{caponeformicagiovanonlanal2013,fioforgogakoparakoNA,anatrielloformicagiovajmaa2017}.

Recall that the Orlicz space $L(N)$ is said to be {\it exponential}
if there exists $\delta>0$ such that the generating Young-Orlicz
function $N = N(u)$ verifies
\begin{eqnarray*}
\lim_{ u \to \infty  } \frac{\ln N(u)}{ [\ln( 2 + u)]^{1+\delta}} =
\infty.
\end{eqnarray*}
For instance, this condition is satisfied when
\begin{eqnarray*}
N(u) = N^{(m)}(u) = \exp \left( |u|^m/m   \right) - 1, \ m = \const
> 0,
\end{eqnarray*}
as well as for an arbitrary Young-Orlicz function which is
equivalent to $N^{(m)}(u)$ or when
$$
N(u) = N_{(\Delta)}(u) \stackrel{def}{=} \exp \left( \  [\ln (1 +
|u|)]^{\Delta} \ \right) - 1, \ \ \Delta = {\const} > 1.
$$

\vspace{4mm}

Denote, as usually, for an arbitrary measurable function $f: X \to
\mathbb R$ its Lebesgue-Riesz norm
$$
||f||_p := \left[ \int_X |f(x)|^p \ d\mu(x)  \right]^{1/p}, \ p \in
[1, \infty).
$$
Suppose that the measure $\mu$ is probabilistic (or, more generally,
bounded): $\mu(X) = 1$. It is known, see e.g. \cite{Ostrovsky Bide
Mart}, that the measurable function $f$ (random variable, \ r.v.)
belongs to the space $L(N^{m}), \ m = {\const} > 0$ iff
$$
\sup_{p \ge 1} \left[ \  ||f||_p \ p^{-1/m} \ \right] < \infty.
$$
Further, the  non-zero function $ f: X \to \mathbb R$ belongs to the
Orlicz space $L(N_{(\Delta)})$ iff, for some non-trivial constant $
C \in (0, \, \infty)$,
$$
\sup_{p \ge 1} \left[ \ ||f||_p \ \exp \left(  - C \ p^{\Delta}
\right)  \ \right] < \infty.
$$

\vspace{5mm}

Define, as usually, for a function $f: X \to \mathbb{R} $ from the
set $M_0(X,\mu)$ its {\it tail function }
\begin{equation}\label{tail function}
 T[f](t) \stackrel{def}{=} \mu \{ \ x: \ |f(x)| \ge t \ \}, \ t \ge 0.
\end{equation}
The function defined in \eqref{tail function} is also known as
\lq\lq distribution function\rq\rq, but we prefer the first name
since the notion \lq\lq distribution function\lq\lq is very used in
other sense in the probability theory.\par

 An arbitrary tail function is left continuous, monotonically non-increasing, takes values in the interval $[0, \mu(X)]$ if $0 < \mu(X) < \infty$
 and in the semi-open interval $[0,\mu(X))$ if $\mu(X) = \infty$.  Besides,
$$
\lim_{t \to \infty} T[f](t) = 0.
$$
The inverse conclusion is also true: such an arbitrary function is
the tail function for a suitable measurable finite a.e. map $ \ f: X
\to \mathbb{R}, \ $ {\it defined on a sufficiently rich measurable
space.}
\par

The set of all tail functions will be denoted by $ \ W: \ $
\begin{equation}
 W =\{ \ T[f](\cdot), \ f \in M_0(X,\mu) \ \}.
\end{equation}

There are many rearrangement invariant function spaces in which the
norm (or quasi-norm) of the function $f(\cdot)$ may be expressed by
means of its tail function $T[f](\cdot)$, for example, the
well-known Lorentz spaces. For the detailed investigation of the
Lorentz spaces we refer the reader, e.g., to \cite{Bennett
Sharpley1988, Lorentz1 1950, Lorentz2 1951, SoriaLor 1998,
SteinWeiss1975}. \par

We introduce here a modification of these spaces. Let $\theta =
\theta(t), \ t \ge 0$, be an arbitrary tail function:  $ \theta \in
W. $
The so-called {\it tail quasi-norm} (or for brevity {\it tail
norm})\, $||f||_{\Tail[\theta]}$  \! of a function $f \in
M_0(X,\mu)$, with respect to the corresponding tail function $
\theta(\cdot),$ is
defined by 
\begin{equation}\label{tail quasi norm}
||f||_{\Tail[\theta]} \stackrel{def}{=} \ \inf \{K > 0 \, : \,
\forall t
> 0 \ \Rightarrow \ T[f](t) \le \theta(t/K) \  \}.
\end{equation}
\vspace{3mm} It is easily seen that this functional satisfies the
following properties:
$$
||f||_{\Tail[\theta]} \ge 0; \ \ \ ||f||_{\Tail[\theta]} = 0 \
\Longleftrightarrow \ f = 0 ;
$$
$$
|| c \ f||_{\Tail[\theta]} = |c| \, ||f||_{\Tail[\theta]}, \ \  c =
{\const} \in \mathbb{R}.
$$

 \vspace{4mm}

Correspondingly, the set of all the functions $f$ belonging to the
set $ M_0(X,\mu)$ and having finite value
 $ ||f||_{\Tail[\theta]} $ is said to be the {\it tail space} $\Tail[\theta].$  \par

The following question is formulated in \cite{Cwikel Kaminska
Maligranda2004} by M. Cwikel, A. Kaminska, L. Maligranda and L.
Pick: {\it \lq\lq  Are the generalized Lorentz spaces really
spaces?\rq\rq}, i.e., can these spaces be normed such that they are
(complete) Banach functional rearrangement invariant spaces? A
particular {\it positive} answer on this question, i.e., under
appropriate simple conditions, may be found in
\cite{OsSirLorNorm2012}. See also \cite[chapter 1, sections
1,2]{OsMono1999}.  \par

\vspace{4mm}

We denote
$$
I(f) = \int f(x) \ d\mu(x)  = \int_X f(x) \ d\mu(x);
$$
if $ \mu$ is a probability measure, we have $ \mu(X) = 1$ and we
replace $(X = \{x\}, \cal{F}, \mu)$ with the standard triplet
$(\Omega = \{ \omega \}, \cal{F}, {\bf P}) $ and, for any numerical-
valued measurable function, i.e., in other words, random variable
 $\xi = \xi(\omega)$, we have
$$
{\bf E} \xi  := I(\xi) = \int_{\Omega} \xi(\omega) \ {\bf P}(d \,
\omega);  \  \ \ T[\xi](t) = {\bf P} ( |\xi| \ge t ), \ t \ge 0.
 $$

\vspace{4mm}

Define now, for an arbitrary Young-Orlicz function $N = N(u)$, the
following tail function from the set $W$

 \begin{equation}
V[N](t) = V_N(t) \stackrel{def}{=} \min \left( \ \mu(X), \ \frac{1}{N(t)} \ \right).
 \end{equation}
Of course, $ \min(c,\infty) = c, \ c \in (0,\infty)$. \par

Suppose $0 \ne f \in  L(N)$; then there exists  a finite positive
constant $ C $ such that $ I (N( |f(\cdot)|/C)) \le 1$; one can take
for instance $C = ||f||_{L(N)}$. \par
It follows from the classical
Markov-Tchebychev's inequality

 \begin{equation}
T[f](t) \le V[N](t/C), \ t \ge 0.
 \end{equation}
In particular,

 \begin{equation}
T[f](t) \le V[N] \left( \ \frac{t}{||f||_{sL(N)}} \ \right), \ t \ge
0.
 \end{equation}

In other words, if $0\ne f \in L(N) < \infty$, then the function
$f(\cdot)$,  as well as its normed version  $\tilde{f} =
f/||f||_{L(N)}$, belongs to the suitable tail space:
\begin{equation}
||f||_{{\Tail}[V[N]]}  \le ||f||_{L(N)} =  ||f||_{sL(N)}.
\end{equation}

\vspace{5mm}

\begin{definition}\label{def weak Orlicz}
 {\rm Let $N$ be a Young-Orlicz function and $f\in M_0(X,\mu)$. We say that $f$ belongs
to the {\it weak} Orlicz space $wL(N)$ and we write $f(\cdot) \in
wL(N)$ iff the following condition is satisfied

\begin{equation}\label{weak orlicz}
||f||_{{\Tail}[V[N]]} < \infty \ \Longleftrightarrow \ f \in
{\Tail}[V[N]].
\end{equation}
We will write for brevity also

$$
||f||_{wL(N)} \stackrel{def}{=} ||f||_{{\Tail}[V[N]]}.
$$
Obviously
\begin{equation} \label{low1}
 ||f||_{wL(N)} \le ||f||_{sL(N)}
 \end{equation}
 }
 and
$$sL(N) \subset wL(N).$$
 \end{definition}

\vspace{2mm}

 \begin{remark}
{\rm Let us emphasize the difference between the general {\it tail
space} $Tail[\theta]$ and the concrete weak Orlicz space $wL(N)$. In
the first case the \lq\lq parameter\rq\rq $\theta$ is an arbitrary
element of the tail set $W$, while for the description of the weak
Orlicz space in the definition \ref{def weak Orlicz} the function
$N(\cdot)$ belongs to the narrow class of Young-Orlicz functions.}
\end{remark}

 The complete review of the theory of these spaces is contained in \cite{Liu MaoFa2013}; see also \cite{Liu Ye2010,Liu Hou MaoFa2017} and the recent paper \cite{Kawasumi2018}.
 It is proved therein, in particular, that these spaces are $F$-spaces and may be normed under appropriate
 conditions, wherein the norm in the corresponding $F$-space or Banach space is linear
and equivalent to the weak Orlicz norm.\par

\vspace{3mm}

{\it There a natural question appears: under what conditions imposed
on the function $N = N(u)$ can the inequality \eqref{low1} be
reversed, of course, up to a multiplicative constant? } \par

\vspace{2mm}

{\bf  In detail, our aim is to find necessary and sufficient
conditions, imposed on the Young-Orlicz function $ N(\cdot)$, under
which}

\begin{equation} \label{defY(N)}
Y(N) \stackrel{def}{=} \sup_{0 \ne f \in \ wL(N) }  \left\{ \
\frac{||f||_{sL(N)}}{||f||_{wL(N)}} \  \right\} < \infty.
\end{equation}

{\bf It is also interesting, by our opinion, to calculate the exact
value of the parameter  $ \ Y(N) \ $ in the case of its finiteness;
we will make this computation in Section 3.} \par

\vspace{4mm}

\begin{remark} { \rm The {\it lower bound}  in the last relation, namely,
$$
\underline{Y}(N) \stackrel{def}{=} \inf_{0 \ne f \in \ wL(N) }
\left\{ \ \frac{||f||_{sL(N)}}{||f||_{wL(N)}} \  \right\},
$$
is known and  $\underline{Y}(N) = 1$.  In detail, it follows from
\eqref{low1} that $\underline{Y}(N) \le 1; $ on the other hand, both
the norms coincide for the arbitrary indicator function of a
measurable set $A$ having a non-trivial measure: $0 < \mu(A) <
\infty$ (see \cite{Liu MaoFa2013}). \par }
\end{remark}

\vspace{2mm}

The comparison theorems between weak as well as between ordinary
(strong) Orlicz spaces and other spaces are obtained, in particular,
in \cite{Bennett Sharpley1988,Buld Mush
OsPuch1992,Iaffei1996,KozOsSir2017,Masta2016,SoriaLor
1998,SteinWeiss1975}, etc.\par

\vspace{2mm}

In both the next examples the space $(X = \{x\}, \cal{F}, \ {\bf
P})$ is probabilistic; one can still assume that $X = [0,1]$,
equipped with the ordinary Lebesgue measure $d\mu(x) = dx$. \par

\vspace{4mm}

\begin{example} {\bf A negative case.}\par
\vspace{4mm}

 {\rm Let $ N(u) = N_p(u) = |u|^p,\ p = \const > 1$; in other words, the Orlicz space $ L(N_p)$ coincides with the
classical Lebesgue-Riesz space $ L_p$:
$$
|\xi|_p = \left[ \ {\bf E}|\xi|^p \ \right]^{1/p}.
$$
The corresponding tail function has the form
$$
 V[N_p](t) = \min(1, \ t^{-p}), \ t > 0.
$$
 On the other hand, let us introduce the r.v. $\eta $ such that
$$
T[\eta](t) := V[N_p](t), \ \  t > 0;
$$
 then, the r.v. $\eta$ has unit norm in the corresponding weak Orlicz space $wL(N_p)$ but

$$
{\bf E}|\eta|^p = {\bf E}\eta^p = p \int_1^{\infty} t^{p-1} \ t^{-p} \ dt =  p \int_1^{\infty} t^{-1} \ dt = \infty,
$$

  $ ||\eta||_p = \infty. \ $ In other words $ Y(N_p) = \infty$.

 As usual, the classical Lebesgue-Riesz norm $||\eta||_p$, \ $ p\ge 1$, of the random variable
 $\eta$ is defined by
$$
 ||\eta||_p \stackrel{def}{=} \left[ {\bf E} |\eta|^p \
\right]^{1/p}.
$$

  }
\end{example}

 \vspace{4mm}

\begin{example} {\bf A positive case.}\par
\vspace{4mm}

{\rm Let now
$$
 N(u) = N^{(2)}(u) = \exp \left(u^2/2 \right) - 1, \ u \in \mathbb{R},
$$
the so-called subgaussian case. It is well-known that the non-zero
r.v. $\zeta$ belongs to the Orlicz space $L(N^{(2)})$ if and only if
there exists $C = {\const} > 0$ such that
$$
T[\zeta](t) \le \exp (- C \ t^2), \ t \ge 0,
$$
or equally
$$
\sup_{p \ge 1} \left[ \ ||\zeta||_p/\sqrt{p} \ \right] < \infty.
$$
Thus, in this case,  $Y(N^{(2)}) < \infty$. \par

The same conclusion holds true also for the  more general so-called
exponential Orlicz spaces, which are in turn equivalent to the Grand
Lebesgue Spaces, see \cite{KozOs1985,KozOsSir2017,OsSirHIAT2007},
\cite[Chapter 1, Section 1.2]{OsMono1999}.

For instance, this condition is satisfied when
\begin{eqnarray*}
N(u) = N^{(m)}(u) = \exp \left( |u|^m   \right) - 1, \ m = {\const},
> 0
\end{eqnarray*}
as well as for an arbitrary Young-Orlicz function which is
equivalent to $N^{(m)}(u)$; or when
$$
N(u) = N_{(\Delta)}(u) \stackrel{def}{=} \exp \left( \  [\ln (1 +
|u|)]^{\Delta} \ \right) - 1, \ \Delta = \const > 1.
$$

}
\end{example}

\vspace{5mm}

\section{Main result.}

\vspace{4mm}

Let $ (X = \{x\}, \cal{F}, \mu)$ be a measurable space with atomless
sigma-finite non-zero measure $ \mu$ and let $N$ be a Young-Orlicz
function. Define an unique value $ t_0= t_0(\mu(X)) \in [0, \infty)$
by
$$
N(t_0) = \frac{1}{\mu(X)};
$$
in particular, when $ \ \mu(X) = \infty, \ $ then $ \ t_0 = 0. \ $  \par

\vspace{3mm}
Denote also
\begin{equation}\label{defJ}
\begin{split}
J(N) &\stackrel{def}{=}  \inf_{C > 0} \int_0^{\infty} N(C \ t) \
\left| \  d V[N](t) \  \right| =\\
&-\sup_{C > 0} \int_0^{\infty} N(C \ t) \ d V[N](t) = -\sup_{C > 0}
\int_0^{\infty} N(C \ t) \ V'[N](t)dt.
\end{split}
\end{equation}

\vspace{2mm}

\noindent Note that the function $t \to V[N](t) $ is monotonically
non-increasing, therefore $|d V[N](t)|=-d V[N](t)$. \par

 \vspace{2mm}
Evidently, when $ t_0 > 0$ we have

$$
\int_0^{\infty} N(C \ t) \ \left| \  d V[N](t) \  \right|
=-\int_{t_0}^{\infty} N(C \ t )\, d \left[ \frac{1}{N(t)} \right].
$$

\vspace{2mm}

\begin{theorem}\label{th equivalence strong-weak}
Let $Y(N)$ and $J(N)$ be defined respectively by \eqref{defY(N)} and
\eqref{defJ}. The necessary and sufficient condition for the
equivalence of the strong and weak Luxemburg-Orlicz's norms, i.e.
$Y(N) < \infty$, is the following:

\begin{equation}\label{condition J(N) finite}
J(N) < \infty,
\end{equation}
or equivalently
\begin{equation} \label{def CN}
\exists C  = C[N] \in (0,\infty) \,: \,  \int_0^{\infty} N(C \ t) \
\left| \  d V[N](t) \  \right| < \infty.
\end{equation}
\end{theorem}

\vspace{4mm}

\begin{remark} \label{even less}

{\rm Evidently, if $ C[N] \in (0,\infty)$, then

$$
\forall C_1 \in (0, C[N]) \ \Rightarrow \int_0^{\infty} N(C_1 \ t) \ \left| \  d V[N](t) \  \right| < \infty.
$$
}
\end{remark}

 \ {\bf Proof.}

 \vspace{3mm}

 \ {\bf A.}  First of all, note that

\begin{equation}\label{note A}
\int_X N(f(x)) \ d\mu(x) \ = - \int_0^{\infty} N(t) \ d\, T[f](t) \
= \int_0^{\infty} N(t) \ |d\, T[f](t)| .
\end{equation}

 \vspace{4mm}

\ {\bf B.  An auxiliary tool.} \par

\vspace{4mm}

\begin{lemma}\label{auxiliary lemma}
Let $ \xi, \ \eta $ be non-negative numerical-valued r.v. such that
$T[\xi](t) \le T[\eta] (t), \ t \ge 0$. Let also $N(u)$ be a
non-negative increasing function, $u \ge 0$. Then
\begin{equation}
{\bf E}N(\xi) \le {\bf E}N(\eta).
\end{equation}
\end{lemma}

\noindent{\bf Proof of Lemma \ref{auxiliary lemma}}.\par
 We can
assume as before, without loss of generality, $X = [0,1]$ with
Lebesgue measure. One can assume also that
$$
\xi(x) = [1 - T[\xi]]^{-1}(x), \  \ \ \eta(x) = [1 -
T[\eta]]^{-1}(x),
$$
where $ \ G^{-1} \ $ denotes a left-inversion for the function $ \
G(\cdot). \ $ Then $ \ \xi(x) \le \eta(x) \ $ and hence $ N(\xi)\le
N(\eta), \ $ and a fortiori $ \ {\bf E} N(\xi) \le {\bf E} N(\eta).
\ $ \par

\vspace{4mm}

 \begin{remark}
 {\rm Of course,  Lemma 2.1 remains true  also for non-finite measure $ \ \mu, \ $  as long as it is sigma-finite.}
  \par
\end{remark}

\vspace{4mm}


 \ {\bf C. \ Necessity.}

 \vspace{4mm}

\noindent Let us introduce the following non-negative
numerical-valued measurable function $g = g(x), \ x \in X$,  for
which
\begin{equation}
T[g](t) = V[N](t), \ t > 0;
\end{equation}
then $ \ g(\cdot)  \in wL(N) \ $ with unit norm in this space. \par

\noindent By the condition $Y(N)<\infty$, the function also $g$
belongs to the space $s L(N)$, therefore

 $$
 \exists C_0  \in (0, \infty) \ \, : \,  \gamma(N) = \gamma_{C_0}(N) \stackrel{def}{=} \int_X N(C_0 \ |g(x)|) \ d\mu(x) <\infty.
 $$
 We deduce, by virtue of \eqref{note A},

 $$
 \int _0^{\infty} N(C_0 \ t) \ |d V[N](t)| = \gamma(N) < \infty,
 $$

 \begin{equation}
J(N) = \inf_{ C > 0} \int_0^{\infty} N(C \ t) \ |d V[N](t)| \le  \int_0^{\infty} N(C_0 \ t) \ |d V[N](t)| = \gamma(N) < \infty.
 \end{equation}

\vspace{5mm}

 \ {\bf D. \ Sufficiency.}\par

\vspace{5mm}

\noindent Assume that the condition $ J(N) < \infty$ is satisfied.
Suppose that the measurable function $ f: X \to \mathbb{R} $ belongs
to the weak Orlicz space $wL(N)$:
\begin{equation}
T[f](t) \le V[N](t/C_2), \ t \ge 0,
\end{equation}
for some finite positive value $C_2$. Let $C_3 = {\const} \in
(0,\infty)$, its exact value  will be clarified below. By using
Lemma \ref{auxiliary lemma} we get

\begin{eqnarray*}
& & \int_X N(C_3 f(x) ) \ d\mu(x) = \int_0^\infty N(C_3t)dT[f](t)
\\
&\leq & \int_0^{\infty} N(C_3\, t) \ |dV[N](t/C_2)| =
\int_0^{\infty} N(C_2 \, C_3 \, t) \ |dV[N](t)| \\
&=& \int_0^{\infty} N(C_4\, t) \ |dV[N](t)| < \infty,
\end{eqnarray*}
if the (positive) value $C_4 := C_2 \ C_3$ is sufficiently small,
for instance $ C_4 \le C[N].$ \par

Thus, the function $ f (\cdot)$ belongs to the strong Orlicz space
$sL(N)$.  $\Box$ \par

\begin{remark}
{\rm The condition of Theorem \ref{th equivalence strong-weak} is
satisfied for the exponential Orlicz space of the form $ \
L(N^{(m)}), \ m > 0, \ $ and is not satisfied for the Orlicz space $
\ L(N_{(\Delta)}), \ \Delta >1, \ $ also exponential space.}

\end{remark}

\vspace{5mm}

\section{Quantitative estimates.}

\vspace{4mm}

It is interest, by our opinion, to obtain the {\it quantitative}
estimation of the constant which appears in the norm inequality for
the embedding $wL(N) \subset sL(N)$; namely, our aim is to compute
the exact value for $Y(N)$, defined in \eqref{defY(N)}.
\par
\vspace{2mm}

 In detail, let $f: X \to \mathbb{R}$ be some function
from the space $ wL(N)$; one can suppose, without loss of
generality,
\begin{equation}
T[f](t) \le V[N](t), \ t \ge 0  \ \Longleftrightarrow \
||f||_{wL(N)} \le 1.
\end{equation}

Assume also that the condition \eqref{condition J(N) finite} is
satisfied, namely $J(N) < \infty$; we want to find the upper
estimate for the value $||f||_{sL(N)}$. \par

\vspace{4mm}

Let us introduce the variable

\begin{equation}\label{y0}
y_0 = y_0(N, \ \mu(X)) := N^{-1}(1/\mu(X)),
\end{equation}
 so that $ \ y_0(N, \infty) = 0, \ \  y_0(N, 1) = N^{-1}(1)$ and define the function

\begin{equation}\label{Q}
Q(k)= Q[N](k):=\int_{y_0}^{\infty} N(y/k) \ \left|  d \frac{1}{N(y)}
\right|, \ \ k \in (1, \ \infty].
\end{equation}
or equally
$$
Q(k) = \int_1^{\infty} N \left( \frac{N^{-1}(w)}{k} \ \right) \
\frac{dw}{w^2}.
$$
Of course $ \ Q(0+) = \infty, \ Q(\infty) = 0$. \par

Denote also
\begin{equation}\label{def K0}
k_0[N] := Q^{-1}(1)  \in [1, \infty).
\end{equation}
Notice that the finiteness of the value $ \ k_0[N] \ $ is quite
equivalent to the condition $J(N)<\infty$ of Theorem \ref{th
equivalence strong-weak}. \par

\vspace{4mm}

\begin{theorem}
Assume that the condition $J(N)<\infty$ is satisfied. Let $k_0[N]$
be defined by \eqref{def K0}. Then
\begin{equation}\label{exactvalue}
||f||_{sL(N)} \le k_0[N] \ ||f||_{wL(N)},
\end{equation}
and the coefficient $k_0[N]$ is here the best possible. Namely,
\begin{equation}
Y(N)=\sup_{0 \ne f \in wL(N)} \left\{ \
\frac{||f||_{sL(N)}}{||f||_{wL(N)}} \ \right\} = k_0[N].
\end{equation}
\end{theorem}

\vspace{4mm}

In other words, $k_0[N]$ is the exact value (attainable) of the
embedding constant in the inclusion $wL(N) \subset sL(N)$.

Moreover, there exists a measurable function $f_0: X \to \mathbb{R}
$, with $||f_0||_{wL(N)}=1$, for which the equality in
\eqref{exactvalue} holds true:
\begin{equation} \label{attain}
||f_0||_{sL(N)} = k_0[N] \ ||f_0||_{wL(N)}.
\end{equation}

Obviously $k_0[N]=+\infty$ when $J(N)=\infty$.

\vspace{4mm}

{\bf Proof.}

First of all, note that the function  $k \to Q(k),  \ k \in (1, \
\infty)$ is continuous, strictly monotonically decreasing and
herewith
$$
Q(\infty) = \lim_{k \to \infty} Q(k) = 0,
$$
by virtue of dominated convergence theorem; as well as

$$
Q(1+) \stackrel{def}{=} \lim_{k \to 1+} Q(k) = \int_{y_0}^{\infty}
N(y) \ \left| d \frac{1}{N(y)}   \right|  =
 \int_{1/\mu(X)}^{\infty}
z \ \left| d \frac{ 1}{z} \right| = \infty,
$$
and the case when $\mu(X) = \infty $  is not excluded.\par Thus, the
value $k_0[N]$ there exists, is unique, positive, and finite: \ $ \
k_0[N] \ \in (1, \infty)$.

Further, assume that the non-zero measurable function $ \ f: \ X \to
\mathbb{R} \ $ belongs to the weak Orlicz space $ \ wL(N); \ $ one
can suppose, without loss of generality, $ ||f||_{wL(N)} = 1$:

\begin{equation}
T[f](t) \le \min \left( \mu(X), \frac{1}{N(t)}   \right) =: T[g](t),
\end{equation}
where $T[g](t)=V[N](t)$.
\par
We deduce, from the definition of the value $k_0[N]$ and using once
again Lemma \ref{auxiliary lemma},

$$
\int_X N \left( \ \frac{f(x)}{k_0[N]}  \ \right)  \ d\mu(x) \le
\int_X N \left( \ \frac{g(x)}{k_0[N]} \right) \ d\mu(x)  =
$$

\begin{equation}
\int_{y_0}^{\infty} N \left( \ \frac{y}{k_0[N]} \ \right)    \left| \ d \frac{1}{N(y)} \  \right| = Q(k_0[N]) =1,
\end{equation}
therefore

\begin{equation}
||f||_{sL(N)} \le k_0[N] = k_0[N] \ ||f||_{wL(N)}.
\end{equation}

So we proved the {\it upper} estimate; the {\it unimprovability} of
ones follows immediately from the relation

\begin{equation}
||g||_{sL(N)} = k_0[N] = k_0[N] \ ||g||_{wL(N)}.
\end{equation}
In detail:
\begin{equation}\label{eq integral}
\int_X N \left( \frac{g(x)}{k_0[N]} \right) \ d\mu(x) = \int_X N
\left( \ \frac{y}{k_0[N]} \ \right) \ \left| \ d V[N](y) \ \right| =
1
\end{equation}
  in accordance with the choice of the magnitude $k_0[N]$. Therefore
$$
||g||_{sL[N]} = k_0[N] \ \hspace{5mm}
$$
and simultaneously $||g||_{wL[N]} = 1$. So, in \eqref{attain} one
can choose $f_0(x):= g(x)$ (attainability).  \par

\vspace{5mm}
%
%
%
%
%

 \begin{example}
 {\rm
Let $(X = \{x\}, \cal{F}, \mu)$ be a {\it probability} space with
atomless sigma-finite measure $\mu(X)=1$. We define the following
Young-Orlicz function, more precisely, the following family of
Young-Orlicz functions
\begin{equation*}
N(u) = N^{(m)}(u) \stackrel{def}{=} \exp \left( \ |u|^m /m\ \right)
- 1, \ \  m = {\const} >1.
\end{equation*}
The case $m = 2$ is known as subgaussian case. The corresponding
tail behavior for non-zero r.v. $ \xi$, having finite weak Orlicz
norm in the space $\left (X, L \left(N^{(m)} \right) \right)$, has
the form
$$
T[\xi](t) \le \exp(- C(m) \ t^{m'}/m'), \ \ t \ge 0, \ \  m'
\stackrel{def}{=} m/(m-1).
$$

\vspace{4mm}

Let us introduce the following modification of the incomplete
beta-function
$$
B_{\gamma}(a,b) \stackrel{def}{=} \int_{\gamma}^1 t^{a -1} \ (1 -
t)^{b-1} \ dt, \ \  \gamma \in (0,1), \  \ a,b  = {\const} \in
\mathbb{R}, \ b
> 0,
$$
and define the variables $ \theta  = \theta(k,m) := k^{-m}, \ k
> 1$, and the function
\begin{equation*}
\begin{split}
G(\alpha) \stackrel{def}{=} B_{1/2}(-1, 1 - \alpha) &= \int_{1/2}^1
t^{-2} \ (1-t)^{-\alpha} \,dt\\
&=\int_0^{1/2}(1 - z)^{-2}\  z^{-\alpha} \,dz,\  \ \  \alpha < 1.
\end{split}
\end{equation*}
With the change of variable $t=1-z$ we have
$$G(\alpha)=\int_0^{1/2}(1-z)^{-2}z^{-\alpha}\,dz.
$$
Using the Taylor series expansion
$$(1-z)^{-2}=\sum_{n=0}^\infty(n+1)z^n, \ \ \ z\in (-1,1)$$
which converges uniformly at least in the closed interval $[0,1/2]$,
we get
\begin{equation*}
G(\alpha)= \sum_{n=0}^\infty(n+1)\int_0^{1/2} z^{n-\alpha}\,dz,
\end{equation*}
which gives
\begin{equation}\label{G series}
G(\alpha) = \sum_{n=0}^{\infty} (n+1)\frac{2^{-n - 1 + \alpha }}{n+1
- \alpha}.
\end{equation}
By \eqref{y0} we obtain
$$
y_0=y_0^{(m)}= N^{-1}(1)=(m \ln 2)^{1/m}
$$
and by \eqref{Q}
\begin{equation*}
\begin{split}
Q_m(k) & = Q[N^{(m)}](k) = \int_{y_0^{(m)}}^{\infty} \left( e^{y^m \
k^{-m}/m} - 1  \right)  \ \left| \ d_y \frac{1}{e^{y^m/m} - 1} \
\right|\\ \\
&=\int_{(m \ln 2)^{1/m}}^{\infty} \left( e^{y^m \
k^{-m}/m} - 1 \right) \frac{e^{y^m/m}\,y^{m-1}}{(e^{y^m/m}-1)^2}\, dy.\\
\end{split}
\end{equation*}
Now we put $x=e^{y^m/m}$, so $dx=e^{y^m/m}\,y^{m-1}\, dy$, \,
$x\in(2,\infty)$; then
\begin{equation*}
\begin{split}
Q_m(k) & =\int_2^\infty(x^{k^{-m}}-1)(x-1)^{-2} \, dx.
\end{split}
\end{equation*}
We make another change of variable $t=1-1/x \ \ \Rightarrow \ \
x=1/(1-t) \ \ \Rightarrow \ \ dx =\frac{dt}{(1-t)^2}$, which yields
\begin{equation*}
\begin{split}
Q_m(k) & =\int_{1/2}^1\frac{1-(1-t)^{k^{-m}}}{(1-t)^{k^{-m}}} \,
\left(\frac{t}{1-t}\right)^{-2}\frac{dt}{(1-t)^2}\\ \\
&=\int_{1/2}^1 t^{-2}\ (1-t)^{-k^{-m}}\, dt-\int_{1/2}^1t^{-2}\,dt=\int_0^{1/2} (1 - z)^{-2} \ z^{-k^{-m}} \, dz - 1\\ \\
&=G(k^{-m})-1=G(\theta(k,m))-1
\end{split}
\end{equation*}
Therefore, the value $k_0=k_0 \left[N^{(m)} \right]=Q^{-1}(1)$
defined in \eqref{def K0} may be found as follows. Define an
absolute constant $ \beta_0$ by means of the relation
\begin{equation}
\int_0^{1/2}  (1 - z)^{-2} \ z^{-\beta_0} \, dz= G(\beta_0) = 2;
\end{equation}
then
$$
\beta_0 \approx 0.431870
$$
and
\begin{equation}
 k_0=k_0 \left[N^{(m)} \right] = \left[\beta_0 \right]^{-1/m},
\end{equation}
 or equally
\begin{equation}
G(k_0^{-m})= 2.
\end{equation}
Note in addition that $ G(0) = 1 $, $G(1^-) = \infty$ and $G$ is
strictly increasing in $(0,1)$, therefore the value $\beta_0$ there
exists and it is unique.
\par
Note that
$$
G(\alpha) > \int_{1/2}^1 (1 - t)^{-\alpha} \ dt = \frac{2^{\alpha -
1}}{1 - \alpha}, \ \alpha \in (0,1),
$$
and, when $ \alpha \to 1^-$,
$$
G(\alpha)  \sim \frac{1}{1 - \alpha}.
$$
If $ \alpha \to 0^+ $, by Taylor series expansion we have
$$
G(\alpha) \sim 1 + C_5  \alpha,
$$
where
$$
C_5 \stackrel{def}{=} \int_0^{1/2} \frac{|\ln z|}{(1 - z)^2} \ dz =
2 \ \ln 2  \approx 1.38629.
$$
Indeed, we put
$$
C_6(\varepsilon) := \int_{\varepsilon}^{1/2} \frac{\ln z}{(1 - z)^2}
\ dz,
$$
so that
$$
C_5 = - \lim_{\varepsilon \to 0+} C_6(\varepsilon), \ \ \varepsilon
\in (0, 1/2).
$$
By means of integration by parts we get
$$
C_6(\varepsilon) =  \int_{\varepsilon}^{1/2} \frac{\ln z}{(1-z)^2} \
dz = \int_{\varepsilon}^{1/2} \ln z \ d\frac{1}{1-z}=
$$
$$
\frac{\ln(1/2)}{1/2} - \frac{\ln \varepsilon}{(1 - \varepsilon)} -
C_7,
$$
where
$$
C_7 = \int_{\varepsilon}^{1/2}\frac{1}{z(1 - z)} \ dz=
\int_{\varepsilon}^{1/2}\frac{dz}{z} + \int_{\varepsilon}^{1/2}
\frac{dz}{(1 - z)} =
$$
$$
\ln(1/2) - \ln \varepsilon - \ln(1/2) + \ln(1 - \varepsilon) = \ln(1
- \varepsilon) - \ln \varepsilon.
$$
Therefore
$$ C_6 = - 2 \ln 2 - \frac{\ln \varepsilon}{1 - \varepsilon} -
\ln(1 - \varepsilon) + \ln \varepsilon \,\to\, - 2 \ln 2,
$$
as long as $ \ \varepsilon \to 0+. \ $ Thus, $ \ C_5 = 2 \ln 2. $

\vspace{5mm}

Note that $\displaystyle\lim_{m \to \infty} k_0 \left[N^{(m)}
\right] = 1$.\par

\vspace{5mm}
To summarize: denote
\begin{equation}\label{inf k0}
\underline{k_0} := \inf_N k_0[N],
\end{equation}
where $ \ " \inf " \ $ in \eqref{inf k0} is calculated over all the
Young-Orlicz functions $N(\cdot)$. We actually proved that
\begin{equation}
\underline{k_0} = 1.
\end{equation}
In detail, it follows from \eqref{low1} that
\begin{equation}
\underline{\tau_0} \ge 1.
\end{equation}
On the other hands,
$$
\underline{\tau_0} \le  \lim_{m \to \infty} \tau_0 \left[N^{(m)}
\right]  = 1.
$$
\vspace{3mm}
Evidently,
$$
\overline{k_0} := \sup_N k_0[N] = \infty.
$$

}
\end{example}

 \vspace{0.5cm} \emph{Acknowledgement.} {\footnotesize The first
author has been partially supported by the Gruppo Nazionale per
l'Analisi Matematica, la Probabilit\`a e le loro Applicazioni
(GNAMPA) of the Istituto Nazionale di Alta Matematica (INdAM) and by
Universit\`a degli Studi di Napoli Parthenope through the project
\lq\lq sostegno alla Ricerca individuale\rq\rq }.\par

 \ The second author is grateful to M. Sgibnev (Novosibirsk, Russia) for sending his interesting articles. \par

\end{document}